\documentclass[12pt]{article}
\usepackage{amsmath,amssymb,amsthm}
\usepackage{enumerate}
\usepackage{hyperref}
\usepackage{graphicx}
\usepackage{url}
\usepackage[mathlines]{lineno}
\usepackage{dsfont} 
\usepackage{color}
\definecolor{red}{rgb}{1,0,0}
\def\red{\color{red}}
\definecolor{blue}{rgb}{0,0,1}

\definecolor{green}{rgb}{0,.6,0}

\numberwithin{figure}{section}   

\newtheorem{thm}{Theorem}[section]
\newtheorem{cor}[thm]{Corollary}
\newtheorem{lem}[thm]{Lemma}
\newtheorem{prop}[thm]{Proposition}

\newtheorem{quest}[thm]{Question}

\theoremstyle{definition}
\newtheorem{rem}[thm]{Remark}

\theoremstyle{definition}
\newtheorem{defn}[thm]{Definition}

\theoremstyle{definition}
\newtheorem{ex}[thm]{Example}

\def\mtx#1{\begin{bmatrix} #1 \end{bmatrix}}

\def\si#1#2{SI^{(#1)}\!\left(#2\right)}
\def\sca#1#2{SC\!A^ {(#1)}\!\left(#2\right)}

\DeclareMathOperator{\rank}{rank}

\newcommand{\Fnn}{F^{n\times n}}

\newcommand{\Fnnnn}{F^{n\x \cdots\x n}}
\newcommand{\Fnmpq}{F^{n_1\x \cdots\x n_d}}
\newcommand{\Fn}{F^{n}}

\newcommand{\bx}{{\bf x}}

\newcommand{\bzero}{{\bf 0}}

\newcommand{\x}{\times}
\newcommand{\bit}{\begin{itemize}}
\newcommand{\eit}{\end{itemize}}
\newcommand{\ben}{\begin{enumerate}}
\newcommand{\een}{\end{enumerate}}
\newcommand{\beq}{\begin{equation}}
\newcommand{\eeq}{\end{equation}}
\newcommand{\bea}{\begin{eqnarray*}}
\newcommand{\eea}{\end{eqnarray*}}
\newcommand{\bean}{\begin{eqnarray}}
\newcommand{\eean}{\end{eqnarray}}
\newcommand{\bpf}{\begin{proof}}
\newcommand{\epf}{\end{proof}\ms}
\newcommand{\bmt}{\begin{bmatrix}}
\newcommand{\emt}{\end{bmatrix}}
\newcommand{\ms}{\medskip}

\newcommand{\beqa}{\begin{array}}
\newcommand{\eeqa}{\end{array}}

\newcommand{\lc}{\left\lceil}
\newcommand{\rc}{\right\rceil}
\newcommand{\lf}{\left\lfloor}
\newcommand{\rf}{\right\rfloor}

\DeclareMathOperator{\nul}{null}

\DeclareMathOperator{\M}{M_0}

\newcommand{\Z}{\operatorname{Z}_0}
\DeclareMathOperator{\Zs}{Z}
\DeclareMathOperator{\Zpd}{Z_{pd}}
\DeclareMathOperator{\I}{I}
\newcommand{\pd}{\gamma_p} 

\newcommand{\s}{\mathcal{S}_{0}} 
\newcommand{\sym}{\mathcal{S}} 
\newcommand{\G}{\mathcal{G}}
\newcommand{\HH}{\mathcal{H}}

\begin{document}

\title{Zero forcing and maximum nullity for hypergraphs}
\author{
   Leslie
Hogben\thanks{Department of Mathematics, Iowa State University,
Ames, IA 50011, USA and American Institute of Mathematics, 600 E. Brokaw Road, San Jose, CA 95112, USA
(hogben@aimath.org).}
 }

\maketitle

 

\maketitle

\begin{abstract}
The concept of zero forcing is extended from graphs to uniform hypergraphs in analogy with the way zero forcing was defined as an upper bound for the maximum nullity of the family of symmetric matrices whose nonzero pattern of entries is described by a given graph:    A family of symmetric hypermatrices is  associated with a uniform hypergraph and zeros are forced in a null vector.  
The value of the hypergraph zero forcing number and maximum nullity are determined for various families of uniform hypergraphs and the effects of several graph operations on the hypergraph zero forcing number  are determined.  The hypergraph zero forcing number is compared to the infection number of a hypergraph and the iteration process in hypergraph power domination. \end{abstract}

\medskip

 \medskip\medskip


The edges of a (simple) graph $G$ describe the nonzero off-diagonal pattern of  symmetric matrices associated with $G$.  The study of the maximum nullity 
of these matrices is an active area of research 
 (see \cite{FH07} and \cite{HLA2ch46}  for surveys).  
 The zero forcing number 
 was introduced in \cite{AIM08} as an upper bound for maximum nullity, and independently in control of quantum systems  \cite{BG07} (precise definitions of these and related terms are given in Section \ref{sMZ}). As noted in \cite{REUF15}, zero forcing is also part of the power domination process on graphs that is used to model the optimal placement of monitoring units in electric networks \cite{HHHH02} (see Section \ref{scompare} for more detail). 

A natural way to extend the idea of maximum nullity from graphs to hypergraphs is to associate a family of symmetric hypermatrices with a uniform hypergraph.   When moving from matrices to hypermatrices there are several possible definitions of rank, nullity, etc. Likewise there are several reasonable choices for how to generalize zero forcing to hypergraphs. 
In Section \ref{sMZ} we present one standard definition of nullity for hypermatrices and the resulting definition of maximum nullity  $\M(H)$ of a hypergraph $H$ that generalizes the definition of zero-diagonal maximum nullity $\M(G)$ of a graph $G$ introduced in \cite{mr0}.  We then  develop a naturally related definition of the zero forcing number $\Z(H)$ of a hypergraph $H$ that generalizes the skew zero forcing number $\Z(G)$\footnote{The skew zero forcing number has traditionally been denoted by $\Zs^-$ but we use $\Z$ to emphasize the connection with zero-diagonal matrices and hypermatrices.} of a graph $G$, extending the graph bound $\M(G)\le \Z(G)$ to hypergraphs, that is, $\M(H)\le \Z(H)$.   In Section  \ref{scompare} we show that $\Z(H)\le \I(H)\le \Zpd(H)$, where $\I(H)$ is the infection number of $H$ introduced in \cite{hyperinfect} and $\Zpd(H)$ is the iterated step in  hypergraph power domination defined by Chang and Roussel in \cite{CR15}.  In Section \ref{sfam} we establish additional properties of  the hypergraph zero forcing number and maximum nullity of hypergraphs, 
including  determining   $\Z(H)$ and $\M(H)$  for several families of hypergraphs and the effects of several graph operations  on the zero forcing number.


\section{Maximum nullity and zero forcing number  for uniform hypergraphs}\label{sMZ}

Throughout $F$ denotes a field.  An 
 {\em order-$d$ hypermatrix} or a $d$-{\em hypermatrix} $A\in\Fnmpq$ is  specified by a $d$-dimensional table of
values and is denoted by 
 $A=[a_{i_{1}\cdots i_{d}}]$ where $a_{i_{1}\cdots i_{d}}$ is the entry in position $(i_1,\dots,i_d)$.  A $d$-hypermatrix $A\in\Fnmpq$ is  {\em hypercubical}  if $n_1=\cdots=n_d$, and this common value is the {\em dimension}.
A hypercubical $d$-hypermatrix $A=[a_{i_{1}\cdots i_{d}}]$ is {\em symmetric} if $a_{i_{\pi(1)}\cdots i_{\pi(d)}}=a_{i_{1}\cdots i_{d}}$ for all $\pi\in S_d$ (where $S_d$ denotes the group of permutations of $\{1,\dots,d\}$).  More information on hypermatrices can be found in \cite{HLA2ch15}.

A  {\em $d$-uniform hypergraph} or a $d$-{\em hypergraph} $H=(V,E)$  has a set of vertices $V$ (also denoted by $V(H)$) and a set of hyperedges $E$ (also denoted by $E(H))$ with each hyperedge being a set of $d$ distinct vertices. 
Since $d$-hypermatrices associate naturally with $d$-hypergraphs, most  hypergraphs discussed here are uniform.
The {\em adjacency matrix} $A(H)\in \Fnnnn$ of a $d$-hypergraph $H$ on $n$ vertices is the symmetric   $d$-hypermatrix that has $a_{i_{1}\cdots i_{d}}=1$ if $\{{i_{1},\dots, i_{d}} \}\in E(H)$ and  $a_{i_{1}\cdots i_{d}}=0$ if $\{{i_{1},\dots, i_{d}} \}\not\in E(H)$ (the latter case includes all subscripts with a repeated index). More information on hypergraphs can be found n \cite{hypergraph-book}.

In this section we define a family of symmetric hypermatrices described by a  hypergraph, the nullity of a hypermatrix, and the maximum nullity of  hypermatrices described by a  hypergraph. 
We define the zero forcing number for a  hypergraph and show that it is an upper bound for maximum nullity.  

\subsection{Hypermatrix nullity and maximum nullity for uniform hypergraphs}\label{ssM}

We begin with a  review of the definitions of sets of symmetric matrices associated with a  graph (all graphs discussed are simple).  Let 
$A=[a_{ij}]\in \Fnn$ be a symmetric matrix. The {\em graph} $\G(A)$ of  $A$ has $V(\G(A))=[n]$ (where $[n]=\{1,\dots,n\}$)   and $\{i,j\}\in E(\G(A))$ if and only if $i\ne j$ and $a_{ij}\ne 0$.   
Let $G=(V,E)$ be a graph. The {\em set of symmetric matrices described by}  $G$ is 
$\sym(G)=\{A\in\Fnn: A \mbox{ is symmetric and }\G(A)=G\}.$
The {\em set of zero-diagonal symmetric matrices described by} a graph $G$ is 
\[\s(G)=\{A\in\Fnn: A \mbox{ is symmetric, $a_{ii}=0$ for $i=1,\dots,n$, and }\G(A)=G\}.\]
Zero-diagonal symmetric matrices described by a graph $G$ can be thought of as weighted adjacency matrices of $G$ and were studied in \cite{mr0}.

A hypercubical $d$-hypermatrix $A=[a_{i_1\cdots i_d}]\in\Fnnnn$ is {\em graphical} if $A$ is symmetric and $a_{i_1\cdots i_d}=0$ whenever ${i_1,\cdots, i_d}$ are not all distinct.  The {\em hypergraph $\HH(A)$ of a graphical $d$-hypermatrix $A$} of dimension $n$ has $V(\HH(A))=[n]$ and $\{i_1,\dots,i_d\}\in E(\HH(A))$ if and only if $a_{i_1\cdots i_d}\ne 0$.  The {\em set of graphical matrices  described by} a $d$-hypergraph $H$ is 
\[\s(H)=\{A\in\Fnnnn: A \mbox{ is graphical and }\HH(A)=H\}.\]

We choose to require the diagonal to be zero in the definition of a graphical hypermatrix for a variety of reasons: It is more natural to have $a_{ii\dots i}=0$ given that $a_{i_1\cdots i_d}=0$ whenever $\{i_1,\dots,i_d\}$ contains any repetition (which is necessary to obtain a uniform hypergraph).  It means that a graphical hypermatrix can be viewed as a weighted adjacency matrix of a hypergraph.   
And it is also related to our definition of null vector below; this  is discussed further after defining a null vector.  The choice to require diagonal elements to be zero means that we are generalizing  $\s(G)$ rather than $\sym(G)$.
 
A vector $\bx\in\Fn$ is a {\em null vector} of a matrix $A\in\Fnn$ if $A\bx=\bzero$, i.e., $\sum_{j=1}^na_{ij}x_j=0$ for every $i=1,\dots,n$.  
There are various possible ways to extend this definition to hypermatrices, and we  choose the next definition.  

\begin{defn}\label{def:nullvec} Let  $A=[a_{i_1\cdots i_d}]\in\Fnnnn$ be a symmetric  $d$-hypermatrix.  A vector $\bx\in\Fn$ is a {\em null vector} of $A$ if 
\beq\label{eq:nullvec} \sum_{j=1}^na_{i_1\cdots i_{d-1}j}x_j=0\ \mbox{ for every submultiset $\{i_1,\dots,i_{d-1}\}\subset [n]$}\eeq   (it is not assumed the values of $i_1,\dots,i_{d-1}$ in \eqref{eq:nullvec} are distinct). The {\em kernel} of $A$, denoted by $\ker A$, is the vector space of null vectors of $A$, and the {\em nullity} of $A$, denoted by $\nul A$, is the dimension of $\ker A$.
The {\em maximum nullity}  of a $d$-uniform hypergraph $H$ is
\[\M(H)=\max\{\nul A : A\in \s(H)\}.\]
\end{defn}

The  definition of null vector could have been stated for a nonsymmetric hypercubical hypermatrix.  However the question then arises as to why the sum is on the last index (the $d$th flattening as defined below).  For a symmetric hypermatrix $A$, $ \sum_{j=1}^na_{i_1\cdots i_{d-1}j}x_j=0$ $ \forall i_1,\dots,i_{d-1}$  is equivalent to  $  \sum_{j=1}^na_{i_1\cdots i_{\ell-1} ji_{\ell+1}\cdots i_{d-1}}x_j=0 \ \forall i_1,\dots,i_{d-1}$, so that question is moot.

This definition of a null vector also provides another reason to restrict to zero-diagonal for $d$-hypermatrices with $d\ge 3$:  Consider a symmetric $d$-hypermatrix $A\in\Fnnnn$ in which $a_{i_1\cdots i_d}=0$ whenever  $2\le |\{i_1,\dots,i_d\}|\le d-1$  but $a_{ii\dots i}\ne 0$ is allowed, aligning with the case of graphs where the diagonal is ignored in the definition of $\sym(G)$.  
Let $\bx\in\ker A$.  Then $a_{kk\dots k}\ne 0$ implies $x_k=0$, so the values of $a_{i_1\cdots i_{d-1}k}$ for $(i_1,\cdots, i_{d-1})\ne (k,\dots ,k)$ are irrelevant when determining the null vector.

This definition of null vector of a symmetric  $d$-hypermatrix $A$  is equivalent to taking a null vector of the transpose of the  $d$th flattening of $A$:  The {\em $d$th flattening} of $A$, denoted by  $\flat_d(A)$, is the $n\x n^{d-1}$ matrix whose $j,i$-entry is $a_{i_1\dots i_{d-1}j}$ where $(i_1,\dots, i_{d-1})$ is the $i$th entry in the lexicographically ordered list of elements in $[n]\x\dots \x[n]$ (with $d-1$ copies of $[n]$); see \cite{HLA2ch15} for more information.  For a symmetric $d$-hypermatrix $A\in\Fnnnn$ and $\bx\in\Fn$,  $\bx$ is a null vector of  $A$ (as just defined) if and only if $\bx$ is a null vector of the matrix  $\flat_d(A)^T$ 
in the usual sense, so $\nul A = \nul \flat_d(A)^T$.  
Thus 
the next result is an immediate consequence of  \cite[Prop 2.2]{AIM08} (which is Proposition \ref{AIM08:prop2.2} for matrices, i.e., 2-hypermatrices).

\begin{prop}\label{AIM08:prop2.2} Let $F$ be a field and let  $A\in\Fnnnn$ be a $d$-hypermatrix.
If $\nul A > k$, then for any set $\alpha\subset [n]$ with $|\alpha|=k$ there is a nonzero vector
$\bx=[x_i]\in\ker A$ such that $x_j=0$ for every $j\in \alpha$. 
\end{prop}

For a $d$-hypergraph $H$ and a matrix   $A=[a_{i_1\cdots {i_d}}]\in\s(H)$, Definition \ref{def:nullvec} can be restated using edges:  If 
$\{i_1,\cdots ,i_{d-1},j\}\notin E(H)$, then $a_{i_1\cdots i_{d-1}j}=0$, so
$a_{i_1\cdots i_{d-1}j}x_j=0$,  Thus  a vector $\bx\in\Fn$ is a  null vector of $A\in\s(H)$ if and only if 
\[
\sum_{\{i_1,\cdots ,i_{d-1},j\}\in E(H)}a_{i_1\cdots i_{d-1}j}x_j=0.
\]
For $H$ a $d$-hypergraph and $A\in\s(H)$, this observation about edges and the symmetry of $A$ suggests using a submatrix of the transpose-flattening of $A$ to test for null vectors.
Let $B$ be an $m\x n$ matrix over $F$.  
For $\alpha\subseteq [m]$ and $\beta\subseteq[n]$, the submatrix of $B$ with rows indexed by $\alpha$ and columns indexed by $\beta$ is denoted by $B[\alpha, \beta]$.
Define $A^\flat=\flat_d(A)^T[\alpha, [n]]$ where $\alpha$ is the set of rows indexed by $\{i_1,\dots,i_{d-1}\}$ such that $i_1<\dots<i_{d-1}$ and there exists an edge $e$ of $H$ containing $\{i_1,\dots,i_{d-1}\}$. 
 Then $\bx$ is a null vector of $A\in\s(H)$ if and only if 
\beq\label{eq:nullvec-e} 
A^\flat\bx=\bzero.
\eeq

\begin{figure}[!ht] 
\begin{center}\scalebox{.5}{\includegraphics{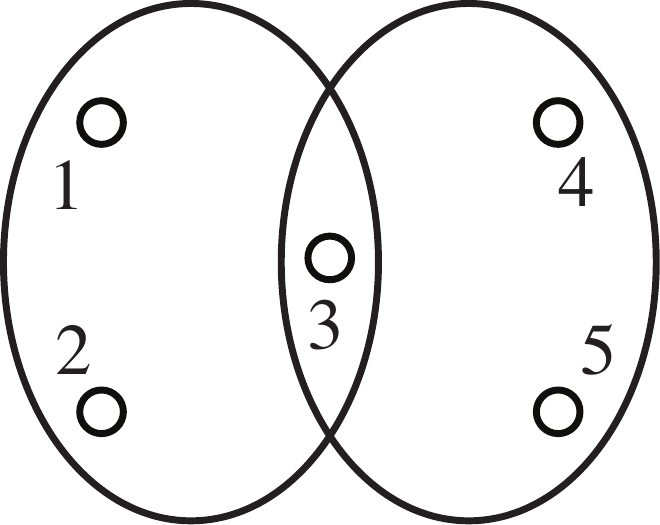}}\qquad\qquad\scalebox{.5}{\includegraphics{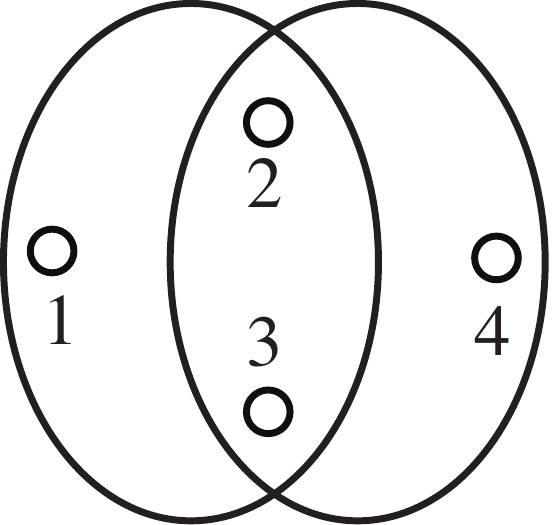}}\\
\ \ \ (a) $H_1$\qquad\qquad\qquad\qquad \ \ (b) $H_2$
 \caption{Hypergraphs $H_1$ and $H_2$ discussed in several examples 
 \label{fig:hyper-ex}}
\end{center}\vspace{-8pt}
\end{figure}

The use of this definition is illustrated in the next two examples.
\begin{ex}\label{ex:M1}  Let $H_1$ be the 3-hypergraph with vertices $\{1,2,3,4,5\}$ and edges $\{\{1,2,3\}$, $\{3,4,5\}\}$ (see Figure \ref{fig:hyper-ex}(a)).  
 For $A=[a_{i_1i_2{i_3}}]\in\s(H_1)$,  $A^\flat=
 \mtx{0 & 0 & a_{123} & 0& 0\\ 0 & a_{123} & 0 & 0& 0\\a_{123} & 0 & 0 & 0& 0\\
0 & 0 & 0 & 0 & a_{345}  \\0& 0& 0 & a_{345} & 0\\ 0& 0 & a_{345} & 0&0}$ with  the rows indexed by $\{1,2\},\{1,3\},\{2,3\},\{3,4\},\{3,5\},\{4,5\}$. Then  $\M(H_1)=0$ because  $\ker A^\flat=\{\bzero\}$.   \end{ex}

\begin{ex}\label{ex:M2}  Let $H_2$ be the 3-hypergraph with vertices $\{1,2,3,4\}$ and edges $\{\{1,2,3\},$ $\{2,3,4\}\}$ (see Figure \ref{fig:hyper-ex}(b)). For $A=[a_{i_1i_2{i_3}}]\in\s(H_2)$,  $A^\flat=\mtx{0 & 0 & a_{123} & 0 \\ 0 & a_{123} & 0 & 0\\a_{123} & 0 & 0& a_{234} \\0 & 0 & a_{234} & 0 \\0 & a_{234} & 0 & 0}$ with  the rows indexed by $\{1,2\},\{1,3\},\{2,3\},\{2,4\},\{3,4\}$. Then  $\M(H_2)=1$ because $[a_{234},0,0,-a_{123}]^T$ is a basis for $\ker A^\flat$. 
 \end{ex}


\subsection{Zero forcing on uniform hypergraphs}\label{ssZ}

We begin with a  review of the definitions of standard and skew zero forcing on a (simple) graph  $G$ as defined in \cite{AIM08} and \cite{IMA10}.  Zero forcing definitions use blue and white vertices, with a blue vertex representing a zero in the null vector of a matrix (in the older literature black is sometimes used instead of blue).  A color change rule allows a white vertex to change color to  blue (but once blue a vertex always remains blue); such a rule is designed by applying information about the matrix and existing zeros in the null vector to conclude that this entry of the null vector must also be zero.  The   {\em standard color change rule}, which is based on  matrices in $\sym(G)$, 
is:  A blue vertex $v\in V(G)$  can change the color of a white vertex $w$ to blue if (1) $ \{v,w\}$ is an edge of $G$,  and  (2) $u$ is white  and  $ \{v,u\}\in E(G)$ implies $u=w$.  The   {\em skew color change rule}, which is based on  matrices in $\s(G)$, 
is:  A vertex $v\in V(G)$  can change the color of a white vertex $w$ to blue if (i) $ \{v,w\}$ is an edge of $G$, and (ii) $u$ is white  and  $ \{v,u\}\in E(G)$ implies $u=w$.  We say {\em $v$ forces $w$} and write $v\to w$  to indicate the color change rule is applied to $w$ by using $v$  to color $w$ blue.
The difference between the two color change rules is that for skew zero forcing, a vertex  need not be blue to force. 

A {\em standard zero forcing set} (respectively, {\em skew zero forcing set}) for $G$ is a set $B\subseteq V$ such that if initially the vertices in $B$ are blue and the vertices in $V\setminus B$ are white, then every vertex can be colored blue by repeated applications of the color change rule. The {\em standard zero forcing number} $\Zs(G)$ (respectively, {\em skew zero forcing number} $\Z(G)$) of $G$ is the minimum cardinality of a standard (skew)  zero forcing set of $G$.    The skew zero forcing number takes its name from its introduction as an upper bound for maximum nullity of skew symmetric matrices described by a graph in \cite{IMA10}, but it also serves as an upper bound for the maximum nullity of zero-diagonal symmetric matrices described by a graph (or any family of matrices with zero-diagonal and off-diagonal nonzero entries described by the edges of the graph) \cite{mr0}.

 Since graphical hypermatrices have all diagonal elements equal to zero, it is the skew color change rule  rather than the standard color change rule that we extend to uniform hypergraphs.  Note that $v$ forcing $w$  in a graph can be interpreted as the elements of edge $\{v,w\}$ other than $w$ forcing $w$; this is the viewpoint we adopt.

\begin{defn}\label{CCR-H} 
Suppose $H$ is a $d$-hypergraph with $d\ge 2$, $B\subseteq V(H)$, every vertex in $B$ is colored blue, and every vertex in $V(H)\setminus B$ is colored white. 
 A set $S\subset V(H)$ of $d-1$ distinct vertices can change the color of white vertex $w$ to blue if\vspace{-5pt}
\ben[(i)]
\item\label{CCR1} $S\cup \{w\}$ is an edge of $H$, and\vspace{-5pt}
\item\label{CCR2} if $u$ is a white vertex and  $S \cup \{u\}$ is an edge of $H$, then $u=w$.\vspace{-5pt}
\een
This is called the  {\em hypergraph color change rule}.  
We say $S$ {\em forces} $w$ and write $S\to w$ to indicate the color change rule is applied to color $w$ blue by  using $S$.
A {\em hypergraph zero forcing set} is a set $B$ such that if  the initial set of blue vertices is $B$, then every vertex can be colored blue by repeated applications of the hypergraph color change rule. The {\em hypergraph zero forcing number} $\Z(H)$ of a hypergraph $H$ is the minimum cardinality of a hypergraph zero forcing set of $H$.    
\end{defn}

We illustrate  hypergraph zero forcing in the next example.

\begin{ex}\label{ex:Z1}  Let $H_1$ be the 3-hypergraph in Example \ref{ex:M1}  (see Figure \ref{fig:hyper-ex}(a)).  Then, $\Z(H_1)=0$ because    $\emptyset$ is a  zero forcing set for $H_1$: $\{1,3\}\to 2$, $\{2,3\}\to 1$, $\{1,2\}\to 3$, $\{3,4\}\to 5$, and $\{3,5\}\to 4$.    \end{ex}

For a hypergraph $H$ and initial set $B$ of blue vertices the {\em derived set of $B$} 
is the set of vertices that are blue after applying the color change rule until no more color changes are possible.

\begin{rem}\label{rem:uniqefinal} As noted in \cite{AIM08} for a graph $G$, the derived set of an initial set $B$  is  unique.  The same reasoning applies to hypergraph zero forcing (and hypergraph infection, and hypergraph power domination zero forcing): 
Any vertex that turns blue
under one sequence of applications of the color change rule can always be turned blue regardless
of the order of color changes:  
Suppose $H$ is a $d$-hypergraph, $B$ is a set of blue vertices, and there is a sequence of forces that results in a derived set $D_1$, and there is another forcing process that colors $D_2$ blue. If $D_1\not\subseteq D_2$, then it is possible to continue forcing: Among vertices in $D_1\setminus D_2$, let $u$ be the first vertex colored blue in the forcing process that produces $D_1$, with $S\to u$.  When $u$ is colored blue, all the vertices that were blue before $u$ are in $D_2$.  Then after the second forcing process produces $D_2$, it is still possible to perform the force $S\to u$ (since coloring additional vertices makes it easier to force).    
Note that the set of forces used to produce  the derived set is usually not unique.
\end{rem}


\begin{thm}\label{thm:MZ} Suppose $H$ is a $d$-hypergraph on $n\ge d\ge 2$ vertices.   
Then \[\M(H)\le\Z(H).\]  \end{thm}
\bpf  We prove the following statement:
\beq\label{eq:zf}  A\in\s(H), \,   \bx\in \ker A,\, B\mbox{ a zero forcing set for } H, \,  \mbox{and}\, x_i=0\,  \forall i\in B \implies \, \bx=\bzero. 
\eeq
Once \eqref{eq:zf} is established, we can choose arbitrarily at most $\Z(H)-1$ zeros in a nonzero vector in $\ker A$ for any $A\in\s(H)$. Then,  $\M(H)\le\Z(H)$ by Proposition \ref{AIM08:prop2.2}.

 Assume that $A\in\s(H)$, $\bx\in\ker A$,   $x_i=0$ for all $i\in B'\subseteq V(H)$, and the color change rule allows the force $S\to w$ with the vertices in $B'$ blue.  Denote the vertices in $S$ by $i_1,\dots,i_{d-1}$.  Then  
 \bea 0&=&\sum_{j=1}^na_{i_1\cdots i_{d-1}j}x_j\\
 &=&\sum_{\{i_1\cdots i_{d-1}j\}\in E(H)}a_{i_1\cdots i_{d-1}j}x_j\\
 &=&\sum_{\{i_1\cdots i_{d-1}j\}\in E(H), j\ne w}(a_{i_1\cdots i_{d-1}j}\cdot 0)+a_{i_1\cdots i_{d-1}w}x_w
 \eea
 Since $\{i_1\cdots i_{d-1}w\}\in E(G)$ implies that $a_{i_1\cdots i_{d-1}w}\ne 0$, necessarily $x_w=0$.
If  $B$ is a zero forcing set for $H$, then by the zero forcing process $\bx=0$. Thus, \eqref{eq:zf} is established. \epf

Theorem \ref{thm:MZ} is applied in the next example.

\begin{ex}\label{ex:Z2}  Let $H_2$ be the 3-hypergraph in Example \ref{ex:M2}  (see Figure \ref{fig:hyper-ex}(b)).  Then, $\Z(H_2)= 1$ because   $\{1\}$ is  a  zero forcing set for $H_2$ with forces $\{1,3\}\to 2$, $\{1,2\}\to 3$, and $\{2,3\}\to 4$, which  implies $\Z(H_2)\le 1$.  Note that $\Z(H_2)\ge 1$ since it was shown that $\M(H_2)=1$ in Example \ref{ex:M1}.  
 \end{ex}

When $H$ is defined as a hypergraph, 
the only color change and zero forcing definitions that apply are those in Definition \ref{CCR-H}, so ``hypergraph" may be omitted from the terminology.  We use the symbol $\Z$ that is associated with skew zero forcing, so for a graph $G$ it does not matter whether we view $G$ as a graph or a 2-hypergraph when writing $\Z(G)$, but for a graph the correct term is ``skew zero forcing.''  The case $d=2$ in Remark \ref{prop:Zub} coincides with the known result $\Z(G)\le n-2$ for a graph $G$ that has an edge \cite{IMA10}.

 \begin{rem}\label{prop:Zub}  Let $n\ge d\ge 2$ and $H$ be a  $d$-hypergraph $H$ on $n$ vertices that has an edge. Then, $\Z(H)\le n-d$, because we can choose any one edge $e=\{w_1,\dots,w_d\}$ and color the remaining $n-d$ vertices blue.  Define $S_i=\{ w_1,\dots,w_{i-1},w_{i+1},\dots, w_d\}$.  Then $S_i \to w_i$ for $i=1,\dots, d$, so $\Z(H)\le n-d$.   
\end{rem}

The {\em degree} of a vertex $v$ of a hypergraph $H$ is the number of edges that contain $v$ and is denoted by $\deg(v)$ (or $\deg_H(v)$ if the hypergraph is not clear).  

 \begin{rem}\label{rem:deg1}  Let $n\ge d\ge 2$ and $H$ be a  $d$-hypergraph $H$ on $n$ vertices.  Suppose $\deg(v)=1$ and let $e=\{v,w_2,\dots,w_{d}\}$ be the edge that contains $v$.  Then every vertex in $e$ except $v$ can be colored blue by the empty set, because $\{v,w_2,\dots,w_{i-1},w_{i+1},\dots,w_d\}\to w_i$.  If edge $e$ has two or more vertices of degree one, then all vertices in $e$ can be colored blue by the empty set.
\end{rem}

Remark \ref{rem:deg1} illustrates a feature of hypergraph zero forcing that is significantly different from skew zero forcing on graphs: In a $d$-hypergraph with $d\ge 3$, a vertex may participate in any number of forces (by combining it with distinct sets of other vertices), whereas in a graph a vertex acts alone to force and thus may perform at most one force.


\section{Comparison of zero forcing and infection and power domination for uniform hypergraphs}\label{scompare} 

In this section we compare the extension of zero forcing to uniform hypergraphs discussed in Section \ref{sMZ} to other extensions of zero forcing to hypergraphs and show our definition is the best upper bound for maximum nullity among these definitions.   
\subsection{Infection  for  hypergraphs}\label{ssI}

Bergen et al.~defined the infection number of a hypergraph as a generalization of the zero forcing number of a graph \cite{hyperinfect}.  In this section we show that the infection number of a uniform hypergraph $H$  is at least as large as the zero forcing number of $H$. 
Suppose $H$ is a  hypergraph 
with  a set $B$ of infected vertices  (and vertices in $V(H)\setminus B$ are uninfected). 
The {\em infection rule} \cite{hyperinfect} allows a non-empty set $S\subseteq B$ of infected vertices to infect all the other vertices in an edge $e\in E$ if\vspace{-4pt}
\ben[(1)]
\item $S \subset e$, and\vspace{-4pt}
\item\label{IF2} if $u$ is an uninfected vertex and $u\not\in e$, then $S \cup \{u\}\not\subseteq e'$\vspace{-4pt}
for every edge $e'$.
\een
An {\em infection set} is a set $B$ such that if initially the set of infected vertices is $B$, then every vertex can be infected by  repeated applications of the infection rule. The {\em infection number} of $H$, denoted by $\I(H)$, is the minimum cardinality of an infection set.

\begin{thm}\label{thm:ZI} Suppose $H$ is a uniform hypergraph. 
Then any infection set for $H$ is a zero forcing set for $H$ and \[\Z(H)\le\I(H).\]  \end{thm}
\bpf  
 Observe that condition \eqref{CCR2} in the hypergraph color change rule could be restated as
\vspace{-4pt}  \ben[(2')]
\item if $u$ is a white vertex and $u\not\in e:= S\cup\{w\}$, then $S \cup \{u\}\not\subseteq e'$
for every edge $e'$ of $H$, \vspace{-4pt} 
\een
which more clearly parallels the infection rule condition \eqref{IF2}.  Differences include that \vspace{-4pt} 
\ben[(a)]
\item a set $S$ of vertices need not be blue/infected to apply the hypergraph color change rule, and \vspace{-4pt} 
\item  a set of maximum cardinality is used for the hypergraph   color change rule.  \vspace{-4pt} 
\een
The first of these properties makes it easier to perform a force, and given that vertices need not be blue to apply the hypergraph color change rule, choosing a maximal set makes it easier to perform a force.  Thus every infection set is a zero forcing set and $\Z(H)\le \I(H)$.
\epf

The next example shows equality is possible in Theorem \ref{thm:ZI}.

\begin{ex}\label{ex:ZeqI}  Let $H$ be the 3-hypergraph shown in Figure \ref{fig:hyper-ex}(b). Then $\I(H)=1$ because  $\I(G)\ge 1$ for every hypergraph $G$ and $\{1\}$ is an infection set:  $\{1\}$ infects \{2,3\} and \{2,3\} infects \{4\}. It was shown in Example \ref{ex:Z2} that $\Z(H)=1$.   \end{ex}

It is also possible to have an arbitrarily large separation between $\Z(H)$ and $\I(H)$  for any $d\ge 3$.   For $d\ge 3$ and $p\ge 2$, define a $d$-hypergraph $S_p^{(d)}=(V,E)$ by 
$V=\{0,1,2,\dots,p(d-1)\}$, and $E=\{e_i:=\{0,(i-1)(d-1)+1,(i-1)(d-1)+2,\dots,i(d-1)\}:i=1,\dots,p\}$.  The hypergraph in Example \ref{ex:Z1} is $S_2^{(3)}$, and  $S_p^{(d)}$ is called a {\em star} because $\displaystyle \left|\cap_{e\in E} e\right|=1$. 

\begin{ex}\label{ex:ZlessI}  Remark \ref{rem:deg1} implies that $\Z(S_p^{(d)})=0$, since every edge has at least two vertices of degree one (because $d\ge 3$).  It is shown that $\I(S_p^{(d)})=p-1$ in  \cite{hyperinfect} 
($S_p^{(d)}$ is one of the graphs 
there called a {flower}).
 \end{ex}

\subsection{Power domination zero forcing for  hypergraphs}\label{ssPD}

For graphs it was shown in \cite{REUF15} that power domination as defined in \cite {HHHH02} can be viewed as a domination step followed by a zero forcing process.  Chang and Roussel extended power domination (and more generally $k$-power domination) to hypergraphs in  \cite{CR15}. 
In this section we identify the hypergraph zero forcing process in \cite{CR15} and show that for uniform hypergraphs the associated zero forcing number is at least as large as the  infection number.  

Vertex $w$ is a {\em neighbor} of vertex $v$ if there is an edge that contains both $v$ and $w$.  A {\em power dominating set} of a hypergraph $H=(V,E)$ is a set $D$ of vertices   that observes all vertices according to the 
 {\em observation rules}  \cite{CR15}, where the first rule is applied once and the second rule is applied repeatedly:
\ben[(1)]
\item\label{PD1} A vertex in $D$ observes itself and  all its neighbors. \vspace{-4pt}
\item\label{PD2} If $v$ is  observed and all the unobserved  neighbors of $v$ are in one edge that contains $v$, then all these unobserved neighbors of $v$ become observed as well.
\vspace{-4pt}
\een
The {\em power domination number} $\pd(H)$ is the minimum cardinality of a power dominating set of $H$.  The second rule can be interpreted as a color change rule, here called the {\em power domination color change rule}:
\bit
\item  If all the white neighbors of a blue vertex $v$  are in one edge that contains $v$, then all these white  neighbors of $v$ change color to blue.
\eit
A {\em (hypergraph) power domination zero forcing set} is a set $B\subseteq V(H)$ such that if initially the set of blue vertices is $B$, then every vertex can be colored blue by repeated applications of the power domination color change rule. The {\em power domination zero forcing number} $\Zpd(H)$ of a hypergraph $H$ is the minimum cardinality of a  power domination zero forcing set of $H$. 

\begin{thm} For a  hypergraph $H$, any power domination zero forcing set is an infection set and $ \I(H)\le \Zpd(H).$
\end{thm}
\bpf The power domination hypergraph color change rule is the same as the infection rule with the restriction that $S$ be a single vertex.  Thus, whenever the power domination hypergraph color change rule can be applied, so can the infection rule.
\epf

\begin{cor} For a  uniform hypergraph $H$,  \[\M(H)\le \Z(H)\le \I(H)\le \Zpd(H).\]
\end{cor}

  The next two examples show that it is possible to have all three parameters $\Z, \I$, and $\Zpd$ equal or distinct. 

\begin{ex}\label{ex:IeqZpd}  Let $H$ be the 3-hypergraph shown in Figure \ref{fig:hyper-ex}(b). It was shown in Examples \ref{ex:Z2} and \ref{ex:ZeqI} that $\Z(H)=\I(H)=1$.  Also,  $\Zpd(H)=1$ because  $\I(G)\le \Zpd(G)$ for every hypergraph $G$ and $\{1\}$ is a power domination zero forcing set:  $1$ colors 2 and 3 blue, and 2 colors 4 blue.   \end{ex}

\begin{ex}\label{ex:IlessZpd}  Let $H$ be the 3-hypergraph shown in Figure \ref{fig:IlessZpd}. We show that $\Z(H)=0, \I(H)=1$, and $\Zpd(H)=2$. 

\begin{figure}[!ht] 
\begin{center}\scalebox{.45}{\includegraphics{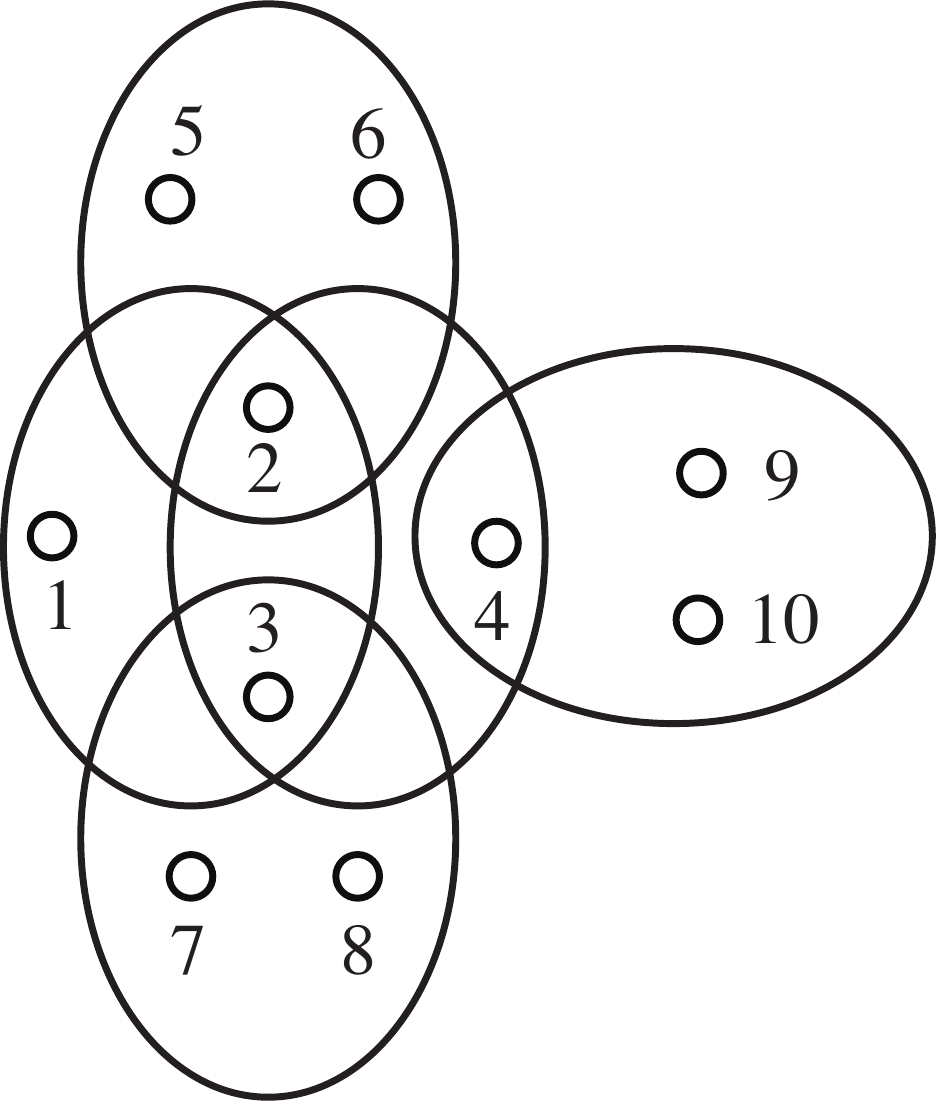}}\\
 \caption{The hypergraph $H$ in Example \ref{ex:IlessZpd}\label{fig:IlessZpd}}
\end{center}\vspace{-8pt}
\end{figure}

\bit
\item $\Z(H)=0$ because edges $\{2,5,6\}, \{3,7,8\}$, and  $\{4, 9,10\}$ each have two degree-one vertices so all vertices except 1 can be colored blue, and then  $\{2,3\}\to 1$.
\item $\I(G)=1$ because $\I(G)\ge 1$ for every hypergraph $G$  and $\{1\}$ is an infection set:  \{1\} infects \{2,3\}, \{2,3\} infects \{4\}, \{2\} infects \{5,6\},  \{3\} infects \{7,8\}, and \{4\} infects \{9,10\}. 
\item $\Zpd(H)=2$ because  we show  no one vertex of $H$ is a power domination hypergraph zero forcing set and \{1,5\} is a power domination hypergraph zero forcing set.
\bit
\item The only one vertex sets we need to consider as possible  power domination hypergraph zero forcing sets are those vertices that in exactly one edge. We consider each in turn.
\bit
\item $D=\{1\}$: 1 colors 2 and 3. 2 has white neighbors in \{2,5,6\} and \{2,3,4\}.   3 has white neighbors in \{3,7,8\} and \{2,3,4\}. 
\item $D=\{5\}$ ($D=\{6\}, \{7\}$ and $\{8\}$ are similar): 5 colors 2 and 6. 2 has white neighbors in \{1,2,3\} and \{2,3,4\}.  6 has no white neighbors.
\item $D=\{9\}$ ($D=\{10\}$ is similar): 9 colors 4 and 10. 4 colors 2 and 3.  2 has white neighbors in \{2,5,6\} and \{2,3,4\}.   3 has white neighbors in \{3,7,8\} and \{2,3,4\}. 10 has no white neighbors.
\eit
\item$D=\{1,5\}$: 1 colors 2 and 3 blue, 5 colors 6 blue, 2 colors 4 blue,  3 colors 7 and 8 blue, and 4 colors 9 and 10 blue.  
\eit
\eit
  \end{ex}
Example \ref{ex:ZlessI} shows that the difference between the infection number and the hypergraph zero forcing number is unbounded for a $d$-hypergraph (independent of $d$), Example \ref{ex:Zpd-I} in  the next section  shows that  the difference between the power domination zero forcing number and the infection number is unbounded if $d$ is allowed to go to infinity. 




\section{Families and further results}\label{sfam}

In this section we study hypergraph zero forcing and maximum nullity further. In particular, we determine the values of these parameters for several families of hypergraphs and study the effect of certain hypergraph operations on the zero forcing number.  We begin with some definitions.

Let $H$ be a $d$-hypergraph.  A {\em $d$-subhypergraph} of $H$  is a $d$-hypergraph $\tilde H$ such that $V(\tilde H)\subseteq V(H)$ and $E(\tilde H)\subseteq E(H)$; in this case $H$ is a {\em $d$-superhypergraph} of $\tilde H$.  A {$d$-subhypergraph} $\tilde H$ of $H$ is {\em induced} if $E(\tilde H)={V(\tilde H)\choose d}\cap E(H)$  where  ${S\choose d}$ denotes the set of  all $d$-element subsets of the set $S$;  $H[U]$ denotes the induced subhypergraph with vertex set $U\subseteq V(H)$.

A {\em path} in a hypergraph $H$ is a vertex-hyperedge alternating sequence \[v_1,e_1,v_2,e_2,\dots,v_s,e_s,v_{s+1}\]  such that $v_1,\dots, v_{s+1}$ are distinct vertices, $e_1,\dots, e_s$ are distinct hyperedges, and $v_i,v_{i+1}\in e_i$ for $i=1,\dots, s$  \cite{hypergraph-book}; such a path is also called a {\em path from $v_1$ to $v_{s+1}$}.
A hypergraph $H$ is {\em connected} if for any two vertices $u$ and $v$ of $H$ there is a path in $H$ from $u$ to $v$.   A {\em connected component} of a hypergraph $H=(V,E)$ is a maximal set of vertices $U$ such that the hypergraph $H[U]$ is connected; in this case, the hypergraph $H[U]$ is called a {\em connected component hypergraph}.  A vertex $v$  is {\em isolated} in a hypergraph $H$ if it is not in any edge of $H$; a connected hypergraph with more than one vertex has no isolated vertices. 

Suppose that the $d$-hypergraph $H$ has connected component hypergraphs $H_1,\dots, H_c$.  It is noted in \cite{hyperinfect} that $\I(H)=\sum_{i=1}^c \I(H_i)$.  Similarly,  it is immediate from the color change rule and the power domination color change rule that  $\Z(H)=\sum_{i=1}^c \Z(H_i)$ and $\Zpd(H)=\sum_{i=1}^c \Zpd(H_i)$.  Let $A$ be a graphical $d$-hypermatrix, and suppose that $\HH(A)$ has  $c$ connected components  $V_1,\dots, V_c$ with connected component hypergraphs $H_i=\HH(A)[V_i], i=1,\dots, c$.  Then $\flat_d(A)^T$ is a block diagonal matrix with the $i$th diagonal block $C_i$ associated with $H_i$, so $\nul A= 
\sum_{i=1}^c \nul C_i$.  Thus, $\M(H)=\sum_{i=1}^c \M(H_i)$. 
Since all of the parameters sum across connected components, it is common to focus on connected  hypergraphs.


\subsection{Zero forcing number and maximum nullity  for families of hypergraphs} 

In this section we study the zero forcing number and maximum nullity  of complete  hypergraphs,  linear  hypergraphs,   interval  hypergraphs (analogous to path graphs),  and circular-arc hypergraphs  (analogous to cycle graphs).
The complete  $d$-hypergraph on $n\ge d$ vertices has all possible edges, i.e.,   $K_n^{(d)}=\left([n], {[n]\choose d}\right)$. 

\begin{prop}\label{prop:Knd-Z}  For $n\ge d$, $\Z(K_n^{(d)})=n-d$.
\end{prop}
\bpf   By Remark \ref{prop:Zub}, $\Z(K_n^{(d)})\le n-d$.
Consider any set $B$ of $n-d-1$ blue vertices, so there are $d+1$ white vertices.  Then any set $S$ of $d-1$ vertices omits at least two white vertices, say $w$ and $u$. Thus $S$ cannot force $w$ because $u$ violates condition \eqref{CCR2}. Therefore, $\Z(K_n^{(d)})\ge n-d$.  \epf

Next we show that $\M(K_n^{(d)})<\Z(K_n^{(d)})$ for $n\ge d+1$.  This parallels the situation for zero-diagonal maximum nullity of complete graphs: $\M(K_n)=n-3<n-2=\Z(K_n)$ \cite{mr0}.
\begin{prop}\label{prop:Knd-M}  For  $n\ge d+1$, $\M(K_n^{(d)})\le n-d-1$. In particular, $\M(K_{d+1}^{(d)})=0$.
\end{prop}
\bpf  
  Define
$r_1= [d+1]\setminus\{d,d+1\}$, $r_2= [d+1]\setminus\{d-1,d+1\}$, $r_3= [d+1]\setminus\{d-1,d\}$,  and $r_{i}= [d+1]\setminus\{d-i+2,d+1\}$ for $i=4,\dots, d+1$. Let $\alpha=\{r_1,\dots,r_{d+1}\}$.  Then 
\[A^\flat[\alpha,[d+1]]=\mtx{
0 & 0 & \dots & 0& 0 & 0 & b_{d+1} & b_d\\
0 & 0 & \dots & 0& 0 & b_{d+1} &0 &  b_{d-1}\\
0 & 0 & \dots & 0& 0  & b_d & b_{d-1}& 0\\
0 & 0 & \dots & 0& b_{d+1}  & 0& 0& b_{d-2}\\
0 & 0 & \dots &  b_{d+1}  &0& 0& 0& b_{d-3}\\
\vdots & \vdots & & \vdots  & \vdots & \vdots & \vdots & \vdots \\
0 & b_{d+1} & \dots &  0  &0& 0& 0& b_2\\
b_{d+1} & 0 & \dots &  0  &0& 0& 0& b_1}\] 
where $b_{i}= a_{1,\dots, {i-1},{i-1},\dots,d+1}$.  From the form of  $A^\flat[\alpha,[d+1]]$ and successive Laplace expansions on the first column, 
\bea \det(A^\flat[\alpha,[d+1]])&=&\pm\, (b_{d+1})^{d-2}\det \mtx{0 & b_{d+1} & b_d\\b_{d+1} &0 &  b_{d-1}\\b_d & b_{d-1}& 0}\\
&=&\pm 2(a_{1,\dots, d})^{d-1}a_{1,\dots, d-2,d-1,d+1}a_{1,\dots, d-2,d,d+1}\ne 0.\eea 
Since  $\rank A^\flat\ge \rank A^\flat[\alpha,[d+1]])= d+1$, $\nul A=\nul A^\flat\le n-d-1$.
\epf

\begin{rem}
 For infection, it is shown in  \cite{hyperinfect} that  $\I(K_n^{(d)})=n-d+1$. 
  The infection set given there (color blue all the vertices outside one edge $e$, and also color one vertex of $e$ blue) also works for power domination zero forcing, and the infection number is a lower bound for power domination zero forcing number, so $\Zpd(K_n^{(d)})=n-d+1$.
 \end{rem}

 A hypergraph is {\em linear} if distinct edges intersect in at most one vertex.  For example, the star $S_p^{(d)}$ is linear.   Every graph is linear when viewed as a 2-hypergraph; however, Corollary \ref{linpath} below shows that the hypergraph zero forcing numbers of linear $d$-hypergraphs with $d\ge 3$  do not behave like those of linear 2-hypergraphs, which are graphs (see \cite{IMA10} for more information about the skew zero forcing numbers of  graphs).  The next, more general, result explains why. 

\begin{prop}\label{smallint} Let $H$ be a 
  $d$-hypergraph  with no isolated vertices such that $|e\cap e'|\le d-2$ for every pair of distinct edges $e$ and $e'$ of $H$.  Then, $\M(H)=\Z(H)=0$. 
\end{prop} 
\bpf
Since $H$ does not have isolated vertices, for each $v\in V(H)$ there is some edge $e_v$ such that $v\in e_v$; let $S_v=e_v\setminus \{v\}$.  Then $|S_v|=d-1> |e_v\cap e'|$ for every  edge $e'\ne e$.  Thus, $S_v\not\subseteq e'$ for $e'\ne e$ and $S_v\to v$.  So $H$ can be forced by the empty set.
That $\M(H)=0$ follows from $0\le \M(H)\le Z(H)$ by Theorem \ref{thm:MZ}.
\epf
 
 Since $|e\cap e'|\le 1$ for a linear $d$-hypergraph, the next result is immediate.
 
\begin{cor}\label{linpath} If $d\ge 3$ and $H$ is a 
linear $d$-hypergraph  with no isolated vertices, then $\M(H)=\Z(H)=0$. 
\end{cor} 

It was established in  \cite{hyperinfect} 
that $\I(H)\le 2$ for a connected linear $d$-hypergraph $H$ in which all vertices have degree at most two   and $\I(H)=1$ if and only if in addition $H$ has a vertex of degree one. For the  zero forcing number of a linear $d$-hypergraph with $d\ge 3$, the degree restriction is unnecessary, and in place of connected we need require only no isolated vertices. However, the restriction on the degree is essential for infection number: The infection number (and therefore also the power domination zero forcing number) is unbounded for linear $d$-hypergraphs as illustrated by  stars, for which $\I(S_p^{(d)})=p-1$   (see Example \ref{ex:ZlessI}).   

Recall that a {\em path graph $P_n$}  is a graph for which there is a linear order of vertices  $v_1,\dots, v_n$ such that $E(P_n)=\{\{v_iv_{i+1}\}: i=1,\dots, n-1\}$ (this is equivalent to defining a path  as already done for  hypergraphs and  defining a path graph to be a graph  whose vertices and edges are those of a path).  Path graphs can be generalized to hypergraphs as interval hypergraphs.

An {\em interval $d$-hypergraph}\,\footnote{The reader is warned that the term {\em interval graph} often means a graph that can be modeled by associating vertices to intervals of the real line with two vertices adjacent if and only if the corresponding intervals overlap.} is a $d$-hypergraph $H$  for which there is a linear order of vertices  $v_1,\dots, v_n$ such that every edge is of the form $\{v_{\ell}, v_{\ell+1}\dots, v_{\ell+d-1}\}$ for some $\ell\in\{1,\dots,n-d+1\}$.  A graph is a connected interval 2-hypergraph if and only if it is a path graph.    It was established in  \cite{hyperinfect} 
that   $\I(H)=1$ for a connected interval $d$-hypergraph $H$.  Thus  the  zero forcing number for an interval $d$-hypergraph is zero or one.  
Example \ref{ex:int1} below provides for each $d\ge 2$ an infinite family of interval $d$-hypergraphs $H$ with  $\M(H)=\Z(H)=1$, and  we show  in Theorem \ref{thm:interval} that all other $d$-hypergraphs have $\M(H)=\Z(H)=0$.

It is convenient to assume that the vertices of an interval hypergraph are denoted by $1,\dots,n$, and we make this assumption.  An interval $d$-hypergraph is determined by its number of vertices $n$ and  its {\em  left endpoint set}  $L(H)=\{\ell_i: e_i=\{{\ell_i}, {\ell_i+1}\dots, {\ell_i+d-1}\}, i=1,\dots,m\}$ where the edges are 
$e_i, i=1,\dots,m$ in order, meaning that $i<j$ implies $\ell_i<\ell_j$.  if $H$ has no isolated vertices, then $\ell_1=1$, and  $\ell_m=n-d+1$.   Since we are interested in generalizing the idea of path, we are interested in connected interval $d$-hypergraphs.  This is equivalent to requiring  that   $\ell_{i+1}\le \ell_i+d-1$.  If $H$ is a connected interval $d$-hypergraph with an edge, then 
$1,e_1,{\ell_2},e_2,\dots,{\ell_m},e_m,n$ is a path that includes every edge of $H$. 
The 3-hypergraphs shown in Figure \ref{fig:hyper-ex} 
are both connected interval 3-hypergraphs on 5 and 4  vertices with left endpoint sets $\{1,3\}$ and $\{1,2\}$, respectively.  Note that for any interval hypergraph on $n$ vertices, $\deg(1)=1=\deg(n)$.


The 3-hypergraph $H_2$ shown in Figure \ref{fig:hyper-ex}(b) and discussed in  Examples \ref{ex:M2} and \ref{ex:Z2} is actually the smallest member of the unique family of interval 3-hypergraphs that have zero forcing number equal to one, as discussed in the next example and theorem.

\begin{figure}[!ht] 
\begin{center}\scalebox{.4}{\includegraphics{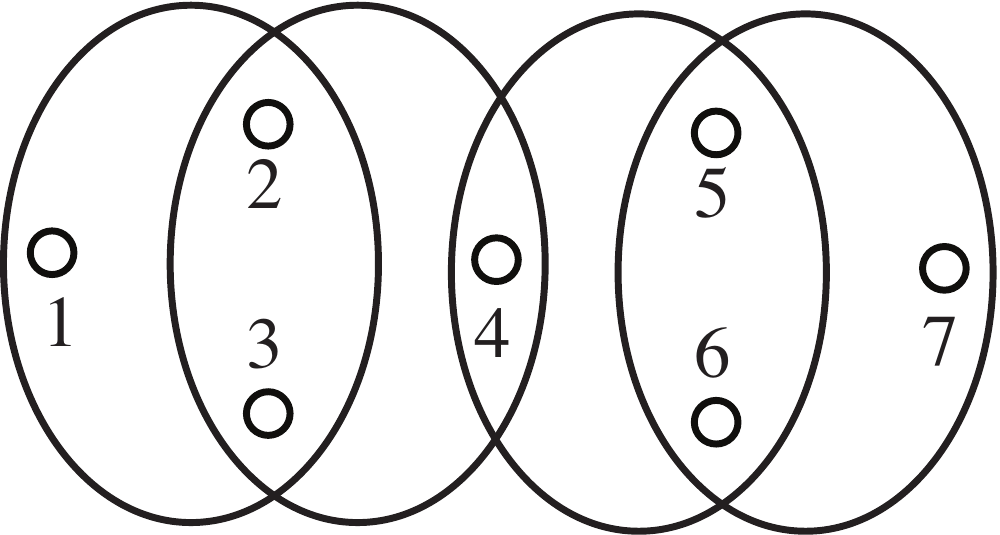}}\\
 \caption{The special interval hypergraph $\si 3 2$ \label{fig:SI32}}
\end{center}\vspace{-10pt}
\end{figure}

\begin{ex}\label{ex:int1} For integers $d\ge 2$ and $s\ge 1$, let $\si d s$ denote the {\em special interval} $d$-hypergraph defined  by $n=sd+1$ and $L(H)=\{(i-1)d+1,(i-1)d+2:i=1,\dots,s\}$; $\si 3 1$ is shown in Figure \ref{fig:SI32}.  Examination of Equation \eqref{eq:nullvec-e} shows that the vector $\bx=[x_j]$ defined by $x_{(i-1)d+1}= (-1)^{i}$ for $i=1,\dots, s$ and all other entries equal to 0 is a null vector of the adjacency matrix of $\si d s$.  Since $\Z(\si d s)\le \I(\si d s)=1$,  $\M(\si d s)=\Z(\si d s)=1$. Note that  $\si d s$ has $2s$ edges and $d$ divides $n-1$. Since  $\si 2 s=P_{2s+1}$, for $d=2$ this is consistent with the known result that maximum nullity and skew zero forcing number of an odd path are one.  \end{ex}

We use a series of lemmas to show that a connected interval $d$-hypergraph  has zero forcing number and maximum nullity equal to zero unless it is a special interval $d$-hypergraph.   

\begin{lem}\label{lem:1zfs} For any connected interval uniform hypergraph $H$ on $n$ vertices, 
 $\{1\}$ and $\{n\}$  are   zero forcing sets for $H$.
\end{lem}
\bpf It is shown in the proof of   \cite[Lemma 4.2]{hyperinfect} 
that $\{1\}$ is an infection set for any such $H$, and a similar argument applies to $n$. The  result then follows from Theorem \ref{thm:ZI}. 
\epf

\begin{lem}\label{lem:interval-emptyforce} Let $H$ be a connected interval $d$-hypergraph determined by $n$ and $L(H)\subset [n]$ with $m$ edges.  For $2\le k \le n-1$, the empty set can force $k$ unless 
 there exists $i$ such that $k=\ell_i=\ell_{i-1}+d-1$ and $\ell_{i+1}=k+1$.  The empty set can force $1$ unless $\ell_{2}=2$, and   
the empty set can force $n$ unless $\ell_{m-1}=n-d$. 
\end{lem}

\bpf
Assume $2\le k \le n-1$.  For $k\in \{\ell_i+1,\dots,\ell_i+d-2\}$ no edge except $e_i$ contains $e_i\setminus\{k\}$.   Thus $e_i\setminus\{k\}\to k$ for $k\in \{\ell_i+1,\dots,\ell+d-2\}$. Thus $k$ can be colored blue by the empty set unless $k$ is in exactly two consecutive edges $e_{i-1}$ and $e_i$ and $k=\ell_i=\ell_{i-1}+d-1$. Furthermore, if $\ell_{i+1}\ne \ell_i+1$, then $\deg (\ell_i+1)=1$ so $e_i\setminus\{\ell_i\}\to \ell_i$.   The cases $k=1$ and  $k=n$ are simpler.
\epf

\begin{lem}\label{lem:interval-dconsec} Let $H$ be a connected interval $d$-hypergraph determined by $n$ and $L(H)\subset [n]$.  
Any set of  $d$ consecutive vertices is a zero forcing set.   
\end{lem}

\bpf  Suppose $B=\{j_1,j_1+1,\dots,j_2-1,j_2\}$ is a maximal set  of consecutive blue vertices with $j_2-j_1\ge d-1$. If  
$j_2=n$, then $B$ is a zero forcing set by Lemma \ref{lem:1zfs}.   So assume 
$j_2+1$ exists;   $j_2+1$ is white by the maximality of $B$. We show $j_2+1$ can be colored blue, and then repeat this procedure until $n$ is colored blue.

If there does not exists an $i$ such that   $j_2+1=\ell_i=\ell_{i-1}+d-1$, then   Lemma  \ref{lem:interval-emptyforce} shows that $j_2+1$ can be colored blue by the empty set.  So assume $j_2+1=\ell_i=\ell_{i-1}+d-1$ for some $i$.  Since $j_2-j_1\ge d-1$, $\ell_{i-1}-1\ge j_1$.  The only edges that can possibly contain  $\{\ell_{i-1},\dots,\ell_{i-1}+d-2\}=e_{i-1}\setminus\{j_2+1\}$ are $e_{i-1}$ and $e_{i-2}$.  If $e_{i-2}$ contains $e_{i-1}\setminus\{j_2+1\}$, then $\ell_{i-2}=\ell_{i-1}-1$, which is blue.  So  $e_{i-1}\setminus\{j_2+1\}\to j_2+1$.  \epf

\begin{thm} \label{thm:interval} Let $d\ge 3$ and let $H$ be a connected interval $d$-hypergraph determined by $n$ and $L(H)\subset [n]$. \[\M(H)=\Z(H)= 
\begin{cases}
1 & \mbox{if $d$ divides $n-1$ and  $H=\si d{\frac{n-1}d}$;}\\
0 & \mbox{otherwise}.
\end{cases}\]
\end{thm}
\bpf  It was shown in Example \ref{ex:int1} that $\M(\si d s)=\Z(\si d s)=1$.   Since $\Z(H)\le \I(H)\le 1$ by Theorem \ref{thm:ZI} and \cite{hyperinfect}, it suffices to show that $\Z(H)\ne 0$ implies $H=\si d s$.
So assume the empty set has failed to color all vertices of $H$ blue and no additional forces are possible.  Then $\ell_1=1$ is white by Lemma \ref{lem:1zfs} and  
$\ell_2=2$ by Lemma \ref{lem:interval-emptyforce}.  Then $e_1\setminus\{k\}\to k$ for $k=2,\dots,d-1$ and $e_2\setminus\{d\}\to d$.  Thus $2,\dots,d$ are blue. Since  there cannot be $d$ consecutive blue vertices by Lemma \ref{lem:interval-dconsec}, $d+1$ is  white, so $\ell_3=d+1$ and $\ell_4=d+2$. Proceeding in order we construct $\si n d$.
  \epf

The conclusion of Theorem \ref{thm:interval} is valid for 2-hypergraphs also, since it is known that $\Z(P_n)=0$ for $n$ even and $\Z(P_n)=1$ for $n$ odd, but the proof given here needs $d\ge 3$.

Interval hypergraphs also behave nicely for power domination zero forcing.
\begin{prop}\label{prop:interval-pd} Let $H$ be a connected interval $d$-hypergraph.  Then $\Zpd(H)= 1$.
\end{prop}
\bpf Note first that $1=\I(H)\le\Zpd(H)$.  Let $B=\{1\}$.  Then $1$ can color blue any white vertices in  $e_1$.    Let $j\ge 2$ and assume that all the vertices in $e_i$ are blue for $i\le j-1$; this implies ${\ell_j}$ is blue.  Observe that ${\ell_j}\not\in e_i$ for $i>j$.  So ${\ell_j}$ can color blue any white vertices of $e_j$. By repeating this process, all the vertices of $H$  can be colored blue. Thus $\Zpd(H)= 1$.
\epf

It is noted in \cite{hyperinfect} that $\I(H)= c$ 
for  an interval $d$-hypergraph $H$ with $c$ connected components, and $\Zpd(H)= c$  by similar reasoning.  The characterization of hypergraph zero forcing number also extends to interval $d$-hypergraphs that are not connected, but we need to consider the types of connected component hypergraphs.

\begin{cor}
Let $H$ be an interval $d$-hypergraph,  
let $c_1$ denote the number of isolated vertices in $H$, and let $c_2$ denote the number of connected components $H_i$ such that $H_i=\si d {s_i}$ for some $s_i$.  Then $\M(H)= \Z(H)=c_1+c_2$. \end{cor}

\begin{figure}[!ht] 
\begin{center}\scalebox{.4}{\includegraphics{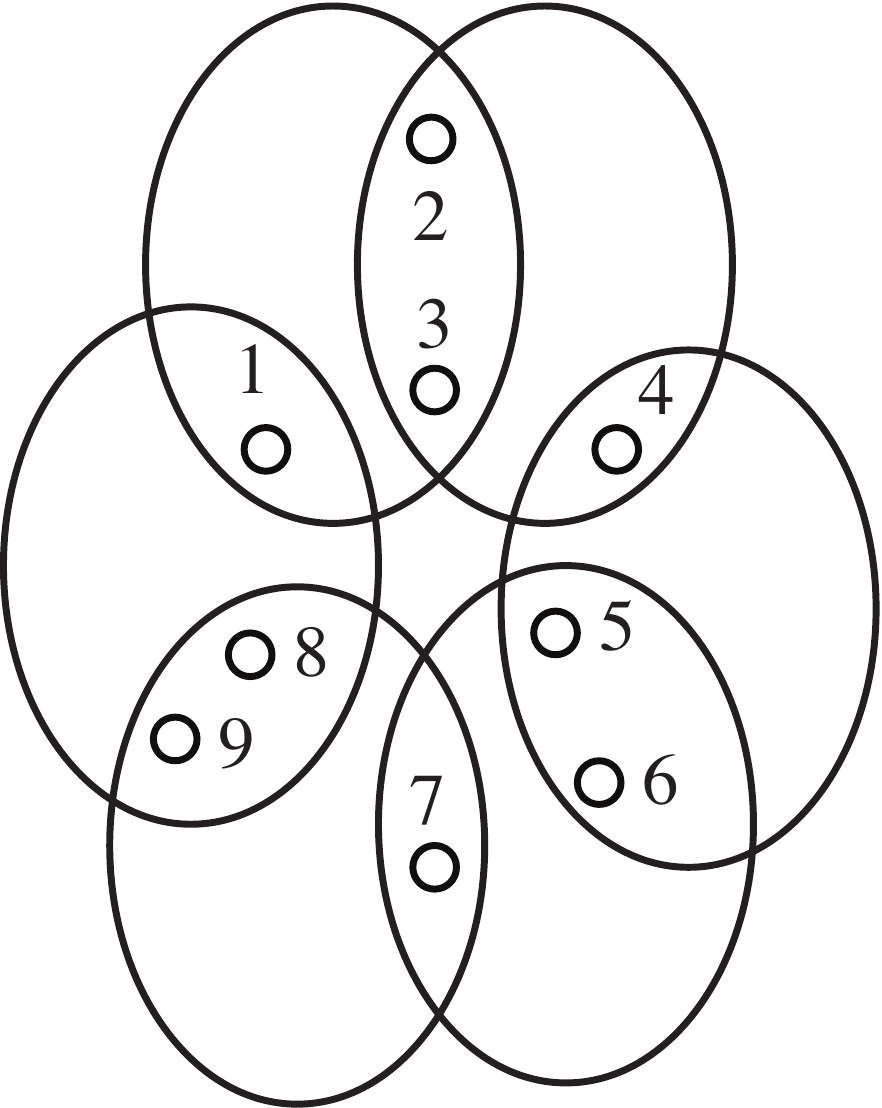}}\\
 \caption{The special circular arc hypergraph $\sca 3 3$ \label{fig:SCA33}}
\end{center}\vspace{-8pt}
\end{figure}

\begin{figure}[!ht] 
\begin{center}\scalebox{.4}{\includegraphics{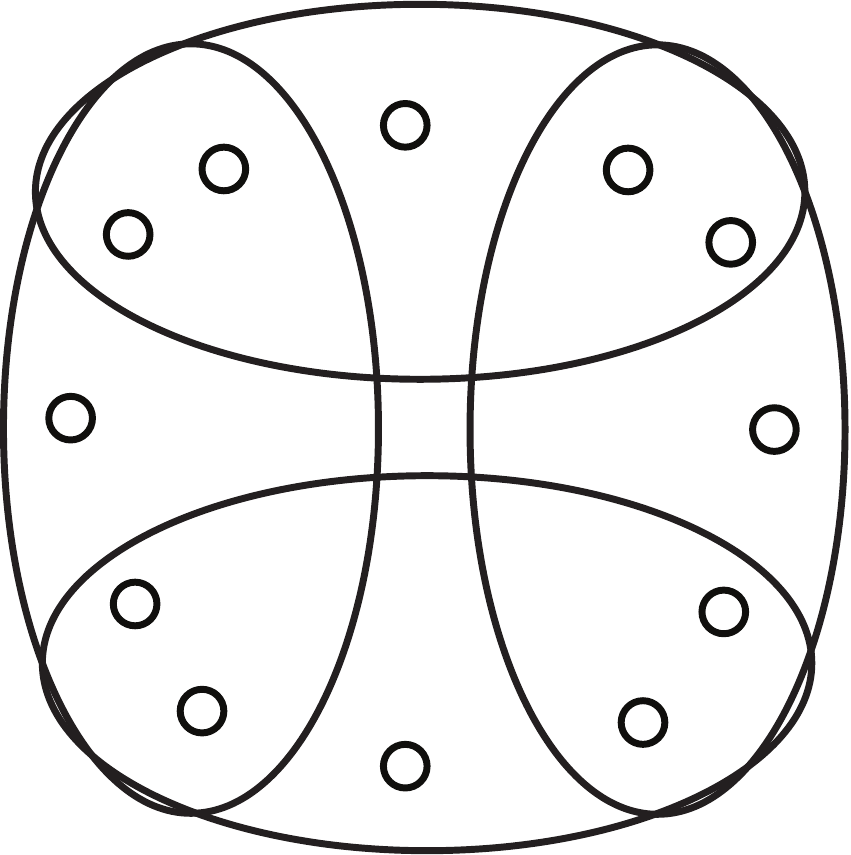}}\\
 \caption{The 2-tight circular arc $5$-hypergraph on 12 vertices 
 \label{fig:Zpd-I}}
\end{center}\vspace{-8pt}
\end{figure}

A {\em cycle} in a hypergraph $H$ is a vertex-hyperedge alternating sequence \[v_1,e_1,v_2,e_2,\dots,v_s,e_s,v_{1}\]  such that $v_1,\dots, v_{s}$ are distinct vertices, $e_1,\dots, e_s$ are distinct hyperedges, and $v_i,v_{i+1}\in e_i$ for $i=1,\dots, s$ (with $s+1$ interpreted as $1$) \cite{hypergraph-book}.
Recall that a {\em cycle graph $C_n$}  is a graph for which there is a cyclic order of vertices  $v_1,\dots, v_n$ such that $E(C_n)=\{\{v_iv_{i+1}\}: i=1,\dots, n\}$ (with $n+1$ interpreted as $1$).    The definitions of cycle and interval hypergraph lead naturally to the idea of a circular-arc hypergraph, which generalizes the idea of cycle graphs to hypergraphs \cite{hyperinfect}.  
A  {\em circular-arc $d$-hypergraph} is a $d$-hypergraph $H$ on $n\ge d+1$ vertices for which there is a cyclic order of vertices  $v_1,\dots, v_n$ such that every edge is of the form $\{v_{\ell}, v_{\ell+1}\dots, v_{\ell+d-1}\}$ for some $\ell\in\{1,\dots,n\}$ (with $n+k$ interpreted as $k$); examples are shown in  Figures \ref {fig:SCA33} and  \ref{fig:Zpd-I}.  

  As with an interval $d$-hypergraph, we assume $V=[n]$, so a circular-arc $d$-hypergraph is determined by $n$ and the {\em first endpoint set} $L(H)=\{\ell_i: e_i=\{{\ell_i}, {\ell_i+1}\dots, {\ell_i+d-1}\}, i=1,\dots,m\}$ where the edges are $e_i,i=1,\dots,m$ and $i<j$ implies $\ell_i<\ell_j$; by convention, $\ell_1=1$.   We are interested in generalizing the idea of a cycle graph and want $1,e_1,{\ell_2},e_2,\dots,{\ell_m},e_m,1$ to be a cycle that includes every edge of $H$.
Thus we also assume  that  $\ell_{i+1}\le \ell_i+d-1$ for $i=1,\dots,m-1$ and $\ell_m\ge n-d+2$.  
For $n\ge d+2$, we obtain a characterization of zero forcing number and maximum nullity that parallels the characterization of these parameters for interval $d$-hypergraphs in Theorem \ref{thm:interval}.

\begin{lem}\label{lem:circ-emptyforce} 
Let $H$ be a circular-arc $d$-hypergraph determined by  $n\ge d+2$ and $L(H)\subset [n]$ with $\ell_{i+1}\le \ell_i+d-1$ for $i=1,\dots,m-1$  and $\ell_m\ge n-d+2$. The empty set can force $k$ unless 
 there exists $i$ such that $k=\ell_i=\ell_{i-1}+d-1$ and $\ell_{i+1}=k+1$.  \end{lem}
 \bpf  Note that for $k\in \{\ell_i+1,\dots,\ell+d-2\}$ no edge except $e_i$ contains $e_i\setminus\{k\}$ because each edge proceeds in cyclic order and the distance going the reverse order is longer since $n\ge d+2$.  Thus $e_i\setminus\{k\}\to k$ for $k\in \{\ell_i+1,\dots,\ell+d-2\}$.  Furthermore, if $\ell_{i+1}\ne \ell_i+1$, then $\deg (\ell_i+1)=1$ so $e_i\setminus\{\ell_i\}\to \ell_i$.   Thus the empty set can force $k$ unless 
 there exists $i$ such that $k=\ell_i=\ell_{i-1}+d-1$ and $\ell_{i+1}=k+1$.  \epf

The proof of the next  lemma is analogous to that of Lemma 
\ref{lem:interval-dconsec}.

\begin{lem}\label{lem:circ-dconsec} Let $H$ be a circular-arc $d$-hypergraph determined by  $n\ge d+2$ and $L(H)\subset [n]$ with $\ell_{i+1}\le \ell_i+d-1$ for $i=1,\dots,m-1$  and $\ell_m\ge n-d+2$.  Any set of  $d$ consecutive vertices is a zero forcing set.   
\end{lem}

For integers $d\ge 2$ and $s\ge 2$, $\sca d s$ denotes the {\em special circular arc} $d$-hypergraph defined  by $n=sd$ and $L(H)=\{(i-1)d+1,(i-1)d+2:i=1,\dots,s\}$; $\sca 3 3$ is shown in Figure \ref{fig:SCA33}. Note that $\sca d s$ can be constructed from $\si d s$ by identifying vertices 1 and $sd+1$.

\begin{thm} \label{thm:circinteval} Let $\ge 3$ and let $H$ be a circular-arc $d$-hypergraph determined by  $n\ge d+2$ and $L(H)\subset [n]$ with $\ell_{i+1}\le \ell_i+d-1$ for $i=1,\dots,m-1$  and $\ell_m\ge n-d+2$.  
\[\M(H)=\Z(H)= 
\begin{cases}
1 & \mbox{if $d$ divides $n$ and  $H=\sca d{\frac n d}$;}\\
0 & \mbox{otherwise}.
\end{cases}\]
\end{thm}

\bpf   
 First consider $H=\sca d{\frac n d}$.  For edge $e=\{i_1,\dots, i_d\}$ let $a_e$ denote the value of $a_{i_1\dots i_d}$ (with entries in all possible orders). Define $A\in \s(\sca d{\frac n d})$ by $a_{e_{2(i-1)+1}}=-1$ and $a_{e_{2i}}=1$ for $i=1,\dots,s$.  Examination of Equation \eqref{eq:nullvec-e} shows that the vector $\bx=[x_j]$ defined by $x_{(i-1)d+1}= 1$ for $i=1,\dots, s$ and all other entries equal to 0 is a null vector of $A$.  Thus $\M(\sca d s)\ge 1$.  Since $d\ge 3$, the empty set can color blue every vertex except those of the form $d(i-1)+1, i=0,\dots,s-1$ by Lemma \ref{lem:circ-emptyforce}.  Let $B=\{1\}$.  Then the $d$ consecutive vertices $1,\dots,d$ are blue, so  all vertices can be colored blue by Lemma \ref{lem:circ-dconsec}.  Thus $\Z(\sca d s)\le 1$, so $\M(\sca d s)=\Z(\sca d s)= 1$.

The proof that $\Z(H)\ne 0$ implies $H=\sca d{\frac n d}$ is analogous to the proof in Theorem  that $H=\si d{\frac {n-1} d}$ if we begin by choosing a white vertex to label as $\ell_1=1$ (after performing all possible forces using the empty set). 
 \epf

The behavior of the zero forcing number of a circular arc $d$-hypergraph graph $H$ with $d\ge 3$ and $n\ge d+2$   differs somewhat from that  of the skew zero forcing number of a cycle, because $\Z(H)\le 1$ and   $\Z(C_n)=2$ for $n$ even (whether viewed as a 2-hypergraph or a graph).

The situation is much more complicated for the infection number, and results in  \cite{hyperinfect} are obtained primarily for a {\em $t$-tight circular-arc $d$ hypergraph}, which is   defined (in our notation) by choosing positive integers $d\ge 3$, $t\le d-1$, and $s\ge \frac{d+1}{d-t} $, and setting $n=s(d-t)$  and $\ell_i=(i-1)(d-t)+1$ for $i=1,\dots,s$; this hypergraph is denoted by $C_n^{(d)}(t)$.   Figure \ref{fig:Zpd-I} shows $C_{12}^{(5)}(2)$. 
Since $\sca d s$ is not tight for $n\ge d+2$,  we know that $\Z(C_{n}^{(d)}(t))=0$ for $n\ge d+2$.  

\begin{prop} \label{tightcircinteval} For positive integers $d\ge 3$ and $t\le d-1$ such that $d-t$ divides $d+1$,  
\[0\le\M(C_{d+1}^{(d)}(t))\le \Z(C_{d+1}^{(d)}(t))= 1.\] 
\end{prop}
\bpf   Note that $\Z(C_{d+1}^{(d)}(t))\le (d+1)-d=1$ by Remark \ref{prop:Zub}. To see that $\emptyset$ is not a zero forcing set, observe that $v_i:=(i-1)(d-t)$ is not in edge $e_i$ for  $i=1,\dots,m:=\frac{d+1}{d-t}$ (with $0$ interpreted as $n$) but is in every other edge, and every other vertex is in every edge. Let $S\subset e_i$ with $|S|=d-1$, so $S$  omits vertex $v_i$ and one other vertex.  If the other vertex is $w\ne v_j$ for $j=1,\dots, m$, then $S$ is a subset only of $e_i$ and $S\to w$.  Thus all vertices  $w\ne v_j$ for $j=1,\dots, m$ can be colored blue.  Now assume $S$ omits both $v_i$ and $v_j$ for some $j\ne i$.  Then $S=e_i\cap e_j$ and each of $e_i$ and $e_j$ has a white vertex not in $S$, so $S$ cannot perform a force.  \epf

The bounds on $\M(C_{d+1}^{(d)}(t))$ are both tight.
If we renumber so that the two vertices that are in both edges are 1 and 3, then Example \ref{ex:M2} is  $C_{4}^{(3)}(1)$, and   $\M(C_{4}^{(3)}(1))=1$. Note that $C_{d+1}^{(d)}(d-1)=K_{d+1}^{(d)}$ so  $\M(C_{d+1}^{(d)})=0$ by Proposition \ref{prop:Knd-M}.  That $\Z(C_{d+1}^{(d)}(t))= 1$ is significantly different from the complicated situation for infection number described in \cite{hyperinfect}, 
where it is shown that $\I(C_{d+1}^{(d)}(d-2))=\frac{d-1}{2}$  when $d$ is odd, implying that the infection number is  unbounded as $d\to \infty$. 
 The next example shows that $\Zpd(H)-\I(H)$ can also be unbounded for tight circular arc $d$-hypergraphs as $d\to\infty$. 

 \begin{ex}\label{ex:Zpd-I} Let $n\ge 2d-1$. It is shown in \cite{hyperinfect} that  $\I(C_{n}^{(d)}(d-1))=2$.   
 We show that $\Zpd(C_{n}^{(d)}(d-1))=d$. Let $B=\{1,\dots, d\}$.  The neighbors of $d$ are $1,\dots,{d-1},{d+1},\dots, {2d-1}$.  Since $1,\dots,{d-1}$ are all blue and ${d+1},\dots, {2d-1}\subset e_2$,  $d$ colors ${d+1},\dots, {2d-1}$ blue.  We can proceed to color the entire hypergraph one edge at a time.  Now suppose $B'\subset V(H)$ contains fewer than $d$ vertices.  If $v\in B'$, then at most $d-2$ other vertices are blue.  Since every vertex has $2d-2$ neighbors, $v$ has $d$ white neighbors, which is too many to include in one edge together with $v$.  So $v$ cannot color any vertex blue.  \end{ex}


\subsection{Hypergraph operations and zero forcing number}

Many  results about the effect of operations on the infection number for hypergraphs \cite{hyperinfect}, zero forcing number for graphs \cite{AIM08, Zspread}, and skew zero forcing number for graphs \cite{fractZF, IMA10} have related versions for the zero forcing number for hypergraphs.  In this section we exhibit several.

Let $H$ be a $d$-hypergraph  with $d\ge 2$.     The $d$-hypergraph $H-e$ is obtained from $H$ by deleting $e$ from the set of edges of $G$.  The $d$-hypergraph $H-v$ is obtained from $H$ by deleting $v$ from the set of vertices of $H$ and deleting from the set of edges every edge that contains $v$. (For $d=2$, these operations are the same as deleting an edge or a vertex from a graph.)  For standard zero forcing on a graph $G$, there are nice bounds relating the zero forcing number of $G$ and the graph resulting from deleting a vertex or an edge: $\Zs(G)-1\le \Zs(G-v)\le \Zs(G)+1$ and $\Zs(G)-1\le \Zs(G-e)\le \Zs(G)+1$ \cite{Zspread}.  Things are not as simple for hypergraph zero forcing with $d\ge 3$.  Since $\Z(S_p^{(d)})=0$  for $d\ge 3$ and deleting the center vertex of $S_p^{(d)}$ leaves $p(d-1)$ isolated vertices, there is no useful upper bound on $\Z(H-v)$ in terms of $\Z(H)$ for $d\ge 3$. 

\begin{prop} Let $H$ be a $d$-hypergraph with $d\ge 2$.  
\ben
\item\label{dele}  $\Z(H)-d\le \Z(H-e)\le \Z(H)+d$.
\item\label{delv} $\Z(H)-1\le \Z(H-v)$.
\item\label{delv2} For $d=2$ (i.e., when $H$ is a graph), $\Z(H-v)\le \Z(H)+1$.
\een
All of these bounds are tight.  \end{prop}
\bpf \eqref{dele}: If $B$ is a zero forcing set for $H$, then $B\cup e$ is a zero forcing set for $H-e$, so $\Z(H-e)\le \Z(H)+d$.  If $\hat B$ is a zero forcing set for $H-e$, then $\hat B\cup e$ is a zero forcing set for $H$, so $\Z(H)\le \Z(H-e)+d$.  A $d$-hypergraph with $d$ vertices and one edge shows the upper bound is tight. A cycle graph on an even number of vertices shows the lower bound is tight for $d=2$, because deleting an edge produces a path graph,  $\Z(C_{2s})=2$, and $\Z(P_{2s})=0$.  

\eqref{delv}: If $\hat B$ is a zero forcing set for $H-v$, then $B'\cup \{v\}$ is a zero forcing set for $H$, so $\Z(H)-1\le \Z(H-v)$.  The complete $d$-hypergraph $K_n^{(d)}$ shows the bound is tight. 

\eqref{delv2}: Let $B$ be a zero forcing set for $H$.  Since $H$ is $2$-hypergraph, $v$ is involved in at most one force. If $\{v\}\to w$, then $B\cup \{w\}$ is a zero forcing set for $H-v$. If $v$ does not perform a force, then $B$  is a zero forcing set for $H-v$.  Thus $\Z(H-v)\le \Z(H)+1$.  A  path graph on an even number of vertices with $v$ adjacent to a vertex of degree one shows the bound is tight. 
\epf

The proof of the next result is adapted from that of Proposition 3.2 in \cite{hyperinfect}, where it is shown that   for any $d$-hypergraph $H$ with  $d\ge 3$, there exists a superhypergraph $H'$ such that $\I(H')=1$.   However, our result applies to 2-hypergraphs also and matches the behavior of the skew zero forcing number on graphs, whereas it is noted in \cite{hyperinfect} that this behavior is very different from the behavior of the zero forcing number on graphs. 

\begin{prop}\label{prop:suphyper0} Let $H=(V,E)$ be a $d$-hypergraph with $d\ge 2$. There exists a $d$-hypergraph $H'$ such that $H=H'[V]$ and $\Z(H')=0$.
\end{prop}
\bpf Let $n$ denote the number of vertices of $H$.  Partition $V$ into $\lf\frac n{d-1}\rf$ sets of $d-1$ vertices, and possibly one additional set of less than $d-1$ vertices; denote these sets by $W_1,\dots, W_k$ where $k=\lc\frac n{d-1}\rc$; if $|W_k|<d-1$ then add  additional vertices to $W_k$ to make its cardinality equal to $d-1$ (it does not  matter which vertices are added).  Define $H'=(V',E')$ by $V'=V\cup \{w_1,\dots, w_k\}$ and $E'=E\cup \{W_i\cup\{w_i\}:   i=1,\dots,k\}$.  Then for $i=1,\dots,k$ and $v\in W_i$, $(W_i\setminus \{v\})\cup\{w_i\}\to v$, and then once every vertex of $V$ is blue, $W_i\to w_i$ for $i=1,\dots, k$.   
\epf

For any set $S$ and object $v$, define $S\x v = \{(x, v):x\in S\} $; $v \x S$ is defined analogously.  The {\em Cartesian product} of $d$-hypergraphs $H$ and $H'$, denoted by $H\Box H'$ is the $d$-hypergraph with
vertex set $V(H) \x V(H')$ and edge set 
\[\{e \x v' : e \in E(H), v' \in V(H')\} \cup \{v \x e' : v\in V(H), e'\in E(H')\}.\] 

\begin{lem}\label{lem:CPforce} Let $H=(V,E)$ and $H'=(V',E')$ be $d$-hypergraphs with $d\ge 3$, $B', S'\subset V'$, $w'\in V'\setminus B'$, and $u\in V$.  Suppose that when the vertices in $B'$ are blue and those in $V'\setminus B'$ are white,  $S'$ can force $w'$ in $H'$.  Then  in $H\Box H'$ when the vertices in $u\x B'$ are blue and those in $u\x(V'\setminus B')$ are white,  $u\x S'$ can force $(u,w')$.
\end{lem}
 \bpf   Since $d\ge 3$, $|S'|\ge 2$.  Thus $u\x S' \not\subseteq e \x v' $ for any $ e \in E$ and $ v' \in V$.  So if  $u\x S'$ is contained in an edge of $H\Box H'$, it is contained in an edge of the form $u\x e'$.  Suppose $u\x S'\cup \{(u,v')\}=u\x e'$ and $(u,v')$ is white. Then $e'=S'\cup \{v'\}$ and  $v'$ is white because all vertices in $u\x  B'$ are blue, so $v\notin B'$. Since $S'\to w'$ in $H'$, this implies $v'=w'$, i.e., $S'\cup \{(u,v')\}=S'\cup \{(u,w')\}$.  Thus $u\x S'\to (u,w')$  in $H\Box H'$. 
\epf

\begin{thm}\label{thm:CP}  Let $H=(V,E)$ and $H'=(V',E')$ be $d$-hypergraphs with $d\ge 3$.  Then 
\[\Z(H\Box H')\le \Z(H)\Z(H').\]
\end{thm}
 \bpf  If $ \Z(H)=0$ or $\Z(H')=0$, then assume without loss of generality that $\Z(H')=0$ and let $\widehat B=\emptyset$.  Otherwise, let $\widehat B=B\x B'$  for  
  minimum zero forcing sets $B$ for $H$ and $B'$  for $H'$.  We show that $\widehat B$ is a zero forcing set for $H\Box H'$.  
 If $\widehat B=\emptyset$, then fix any vertex  $u\in V(H)$.  Otherwise, fix any vertex $u\in B$.  Lemma  \ref{lem:CPforce} implies that
 the zero forcing process on $H'$ can proceed in $H\Box H'$ to color all the vertices in $u\x V(H')$.  Thus  all vertices are blue  if $\widehat B=\emptyset$, and otherwise all vertices in $u\x V'$ are blue  for every $u\in B$.  In the latter case then apply the same reasoning using  the zero forcing process in $H$ to color all the remaining vertices.  
Therefore, $\Z(H\Box H')\le |B| |B'| = \Z(H)\Z(H')$.  
\epf
 
In contrast to the simple upper bound on the zero forcing number of a Cartesian product in Theorem \ref{thm:CP}, the analogous bound for infection number requires defining the 2-infection number, in which at least two infected vertices are needed to infect (see \cite{hyperinfect}). 
Note also that for graphs $G$ (i.e., 2-hypergraphs), it is not generally true that $\Z(G\Box G')\le \Z(G)\Z(G')$:  For example, $\Z(K_2)=0$ but $K_2\Box K'_2=C_4$ and $\Z(C_4)=2$.

For  $d$-hypergraphs $H$ and $H'$ with $d\ge 3$, 
it is immediate that   $\Z(H)=0$ or $\Z(H')=0$ implies $\Z(H\Box H')=\Z(H)\Z(H')=0$, and we do not know of any examples where $\Z(H\Box H')<\Z(H)\Z(H')$.
This naturally raises the next question.  

\begin{quest}\label{q:cartprod}
Is it true that $\Z(H\Box H')= \Z(H)\Z(H')$ for all $d$-hypergraphs $H$ and $H'$ with $d\ge 3$?
\end{quest}

We answer Question \ref{q:cartprod} in the affirmative for hypergraphs with  $\Z(H')= 1$ in the next result, which of course also applies when $\Z(H)= 1$.

\begin{thm}\label{prop:cartprod1}  Let $H=(V,E)$ and $H'=(V',E')$ be $d$-hypergraphs with $d\ge 3$ such that $\Z(H')=1$.  Then 
$\Z(H\Box H')=\Z(H).$
\end{thm}
 \bpf  Let $\tilde B\subset V(H\Box H')$ with $|\tilde B|<\Z(H)$.  Since $\Z(H\Box H')\le \Z(H)$ by Theorem \ref{thm:CP}, it suffices to show that $\tilde B$ cannot color all vertices of  $H\Box H'$ blue.  
Let $B=\{ u\in V(H): \exists u'\in V(H') \mbox{ such that }(u,u')\in\tilde B\}$.  Note that $B\le |\tilde B|<\Z(H)$, so $B$ is not a zero forcing set for $H$.  Let $D$ denote the derived set of  $B$ (in $H$), and let $D'$ denote the derived set of  the empty set   in $H'$. Observe that both $V\setminus D$ and $V'\setminus D'$ are nonempty.  We show that none of the vertices in $(V\setminus D)\x (V'\setminus D')$ are in the derived set of $\tilde B$ (in $H\Box H'$). 

By definition, every edge of $H\Box H'$ is of the form $v\x e'$ or $e\x v'$ where $v\in V, e\in E, v'\in V',$ and $ e'\in E'$.  So if $\tilde S\subset V\x V'$ is a set of $d-1$ vertices that can force $(w,w')$, then $\tilde S\subset v\x e'$ or $\tilde S\subset e\x v'$, but not both because  $d\ge 3$ implies $|S'|\ge 2$.  
By Remark \ref{rem:uniqefinal}, the  order in which we perform the forces does not matter, provided at the end we ensure all options have been tested. So first perform all possible  forces using $\tilde S\subset v\x e'$.  This will color   $(v,u')$ blue  for all $v\in V$ and $u'\in D'$.  Also some additional vertices of the form $(u,v')$ with $u\in B$ and $v\in V'$ will be colored blue.   Assuming additional vertices have been colored blue can not reduce the derived set of $\tilde B$.  So assume  $(u,v')$ is blue for all $u\in B$ and $v'\in V'$.
Observe that every vertex of the form $(y,x')$ with $y\in V\setminus B$ and $x'\in V'\setminus D'$ is still white; this includes every vertex of the form $(x,x')$ with $x\in V\setminus D$ and $x'\in V'\setminus D'$. 
Now perform all possible  forces using $\tilde S\subset e\x v'$.  After this set of forces, every vertex of the form $u\x v'$ is blue for $u\in D$ and $v'\in V'$, but no additional vertices of the form $(x,v')$ for $x\in V\setminus D$ have been colored blue.  Thus $(x,x')$ is still white for every $x\in D$ and $x'\in D'$.  Since no additional vertices in  $x\x V'$ have been colored blue, and all vertices in $u\x V'$ are blue for $u\in D$, no further forces using $\tilde S\subset v\x e'$  are possible, and thus no additional forces are possible.  
\epf

 \begin{cor}\label{cor:cartprod}  Let $H=(V,E)$ and $H'=(V',E')$ be $d$-hypergraphs with $d\ge 3$ such that 
 $\Z(H)\Z(H')\le 3$.  Then 
$\Z(H\Box H')=\Z(H)\Z(H').$
\end{cor}
\bpf Without loss of generality, $\Z(H')\le \Z(H)$.  Since $\Z(H)\Z(H')\le 3$, $\Z(H')\le 1$ and the result follows from Theorems \ref{thm:CP} and \ref{prop:cartprod1}.  \epf


\section{Concluding remarks}

It is not surprising that the zero forcing number defined in Section \ref{sMZ} is a better bound for maximum nullity than the infection number defined in \cite{hyperinfect}, since the authors of \cite{hyperinfect} are explicit that they are not trying to generalize zero forcing as an upper bound for maximum nullity.  They say, ``For hypergraphs, there is no matrix analogous to the adjacency matrix of a graph and this notion of a set of entries in
a vector forcing the other entries to be zero in a proposed null vector does not apply."  While it is true that for a hypergraph there is no {\em matrix} analogous to the adjacency matrix, we would argue that for a uniform hypergraph the natural analog is the adjacency hypermatrix.  
When using hypermatrices,  the issue is not that there are no analogous concepts for nullity and forcing zeros in a null vector, but rather there as several possible choices for nullity of a hypermatrix, which lead naturally to associated concepts of zero forcing.  Here we have made one standard choice for the  definition of nullity that led  to the concept of zero forcing that has been studied.  Other choices for nullity of a hypermatrix will lead to other definitions of zero forcing in hypergraphs, and these are also worthy of study.  

We have obtained a variety of results for the zero forcing number and maximum nullity as defined here, but numerous interesting questions remain.  Here we highlight a few, in addition to Question \ref{q:cartprod}.  There are many additional interesting families of hypergraphs for which the zero forcing number could be determined, and additional hypergraph operations whose effect on zero forcing number could be studied.  
There are also many interesting questions to be studied regarding the power domination zero forcing number, including many of the types of results obtained in \cite{hyperinfect} for the infection number and here for the  zero forcing number, but also the question of a relationship between the hypergraph power domination number defined in \cite{CR15} and power domination zero forcing define here, paralleling that between graph power domination number and zero forcing number \cite{REUF15}.

\section*{Acknowledgments}

I wish to express my deep gratitude  to Edinah Gnang for suggesting that a hypergraph zero forcing number could be defined as an upper bound for hypermatrix nullity, and for many helpful discussions about hypermatrices.   


\end{document}